\documentclass{hha}
\usepackage{amssymb,amsmath}
\usepackage{graphics, xy}

\def\tr{\triangleright}
\def\tl{\triangleleft}

\newcounter{dnum}
\newcounter{tnum}
\newcounter{en}
\newenvironment{theorem}[1][Theorem]{\medskip\noindent{\bf #1 \arabic{tnum}.}\quad}%
{\par\medskip \addtocounter{tnum}{1}}
\newenvironment{definition}[1][Definition]{\medskip\noindent{\bf #1 \arabic{dnum}.}\quad}%
{\par\medskip \addtocounter{dnum}{1}}
\newenvironment{lemma}[1][Lemma]{\medskip\noindent{\bf #1 \arabic{tnum}.}\quad}%
{\par\medskip \addtocounter{tnum}{1}}
\newenvironment{corollary}[1][Corollary]{\medskip\noindent{\bf #1 \arabic{tnum}.}\quad}%
{\par\medskip \addtocounter{tnum}{1}}
\newenvironment{example}[1][Example]{\medskip\noindent{\bf #1 \arabic{en}.}\quad}%
{\par\medskip \addtocounter{en}{1}}

\begin{document}

\title{Matrices and Finite Quandles}
 
\author{Benita Ho}
\email{benitakho@yahoo.com}
\address{University of California, Riverside \\ 900 University Avenue \\ 
Riverside, CA, 92521}
\author{Sam Nelson} 
\email{knots@esotericka.org}
\address{University of California, Riverside \\ 900 University Avenue \\ 
Riverside, CA, 92521 }


\classification{57M27.}

\keywords{Finite quandles.}

\begin{abstract}
\noindent Finite quandles with $n$ elements can be represented as $n\times n$ 
matrices. We show how to use these matrices to distinguish all isomorphism 
classes of finite quandles for a given cardinality $n$, as well as how to 
compute the automorphism group of each finite quandle. As an application, we 
classify finite quandles with up to $5$ elements and compute the automorphism 
group for each quandle.  
\end{abstract}

\received{May 29, 2005}   
\revised{July 9, 2005}    
\published{November 7, 2005}  
\submitted{Ronald Brown}  

\volumeyear{2005} 
\volumenumber{7}  
\issuenumber{1}   

\startpage{197}     

\maketitle

\section{Introduction}

A \textit{quandle} is a set $Q$ with a binary operation $\tr:Q\times Q\to Q$
satisfying the three axioms 
\newcounter{q}
\begin{list}{(\roman{q})}{\usecounter{q}}
\item for every $a\in Q$, we have $a\tr a= a$,
\item for every pair $a,b\in Q$ there is a unique $c\in Q$ such that
$a=c\tr b$, and
\item for every $a,b,c\in Q$, we have $(a\tr b) \tr c = (a\tr c) \tr (b\tr c)$.
\end{list}

The uniqueness in axiom (ii) implies that the map $f_b:Q\to Q$ defined
by $f_b(a)=a\tr b$ is a bijection; the inverse map $f_b^{-1}$ then defines
the \textit{dual operation} $a\tl b= f^{-1}_b (a)$. The set $Q$ then forms a
quandle under $\tl$, called the \textit{dual} of $(Q,\tr)$.

Quandle theory may be thought of as analogous to group theory. Indeed,
groups are quandles with the quandle operation given by $n$-fold 
conjugation for an integer $n$, i.e., \[a\tr b = b^{-n}ab^n.\] Another 
important example of a type of quandle structure is the category of 
\textit{Alexander quandles}, i.e., modules $M$ over the ring 
$\Lambda=\mathbb{Z}[t^{\pm 1}]$ of Laurent polynomials in one variable with 
quandle operation \[a\tr b = ta + (1-t)b.\] The second author has written
elsewhere on Alexander quandles; see \cite{N1} and \cite{N2}.

Other examples of quandles include \textit{Dehn quandles}, i.e., the set of 
isotopy classes of simple closed curves on a surface $\Sigma$ with action 
given by Dehn twists, and \textit{Coxeter quandles}, i. e., 
$\mathbb{R}^n \setminus 0$ with \[ u\tr v = 2\frac{(u,v)}{(v,v)} -u\]
where $(,)$ is a symmetric bilinear form. See \cite{Y} and \cite{FR} for more.

So far, quandles have been of interest primarily to knot theorists, due to 
their utility in defining invariants of knots. In \cite{J}, a quandle is 
associated to every topological space, called the \textit{fundamental 
quandle.} In particular, it is shown that isomorphisms of the \textit{knot 
quandle} (definable from a knot diagram by a Wirtinger-style presentation)
preserve peripheral structure, making the knot quandle a complete
invariant of knot type considered up to homeomorphism of topological pairs, 
though not up to ambient isotopy. 

Finite quandles have been used to define invariants of both knots and links
in $S^3$ and generalizations of knots such as knotted surfaces in 
$\mathbb{R}^4$ and virtual knots. The simplest example of such an 
invariant is the number of homomorphisms from the knot quandle to a
chosen finite quandle. One can also obtain knot invariants by counting 
homomorphisms with crossings weighted by quandle cocyles arising in various 
quandle cohomology theories. See \cite{C} and \cite{CM} for more.

In this paper, we show how to associate to any finite quandle $Q=\{x_1, x_2, 
\dots, x_n\}$ an $n\times n$ matrix $M_Q$. We then define an equivalence 
relation on the set of $n\times n$ quandle matrices, which we call 
$p$-equivalence (after we decided that our initial choice of 
``$\rho$-equivalence'' would be confusing). Our main theorem then says that 
two such matrices represent isomorphic quandles iff they are $p$-equivalent.

We then give an algorithm for applying this result to determine all 
isomorphism classes of quandles with $n$ elements as well as their 
automorphism groups, and as an application we determine all isomorphism 
classes and automorphism groups of quandles with up to $5$ elements.

After the initial posting of the preprint of this paper to arXiv.org, a paper
with similar results was posted by Lopes and Roseman \cite{LR}. We have also 
learned that related results were obtained by Hayley Ryder in her dissertation
\cite{R}.

\section{The matrix of a finite quandle}

\addtocounter{tnum}{1}
\addtocounter{en}{1}
\addtocounter{dnum}{1}

Let $Q=\{x_1, x_2, \dots, x_n\}$ be a finite quandle with $n$ elements. We 
define the matrix of $Q$, denoted $M_Q$, to be the matrix whose entry in 
row $i$ column $j$ is $x_i\tr x_j$:
\[ M_Q = \left[ \begin{array}{cccc}
x_1\tr x_1 & x_1\tr x_2 & \dots & x_1\tr x_n \\
x_2\tr x_1 & x_2\tr x_2 & \dots & x_2\tr x_n \\
\vdots & \vdots & \ddots & \vdots \\
x_n\tr x_1 & x_n\tr x_2 & \dots & x_n\tr x_n \\
\end{array}\right].\]
The matrix $M_Q$ is really just the quandle operation table considered as
a matrix, with the columns acting on the rows. In particular, if the
elements of the quandle are the numbers $Q=\{1,2,\dots, n\}$ with
$M_Q=(\alpha_{ij})$, call $M_Q$ an \textit{integral quandle matrix}. We may
obtain an integral quandle matrix by suppressing the ``$x$''s in the notation 
and just writing the subscripts; hence we lose no generality by restricting
our attention to integral quandle matrices. If the entries on the diagonal 
in an integral quandle matrix are in the usual order, i.e., $\alpha_{ii}=i$, 
then $i\tr j$ is just the entry in row $i$ column $j$. An integral quandle 
matrix of this type is in \textit{standard form}.

The quandle axioms place certain restrictions on what kind of matrices can 
arise from a quandle. 

\begin{lemma} \label{mqdef}
Let $M=(\alpha_{ij})$ be an $n \times n$ matrix with 
$\alpha_{ij}\in \{1,2,\dots,n\}$. 
Then $M=M_Q$ for a finite quandle $Q$ if and only if the following 
conditions are satisfied:
\newcounter{mq}
\begin{list}{(\roman{mq})}{\usecounter{mq}}
\item{The diagonal entries are distinct, i.e. $\alpha_{ii}=\alpha_{jj}$ 
implies $i=j$. If this condition is satisfied, denote the row number 
containing $x$ on the diagonal by $r(x)$ and the column number containing 
$x$ on the diagonal by $c(x)$.}
\item{The entries in each column are distinct, i.e. $\alpha_{ij}=\alpha_{kj}$ 
implies $i=k$.}
\item{The entries must satisfy 
\[\alpha_{r(\alpha_{r(x)c(y)})c(z)} = 
\alpha_{r(\alpha_{r(x)c(z)})c(\alpha_{r(y)c(z)})},\]
or if we denote $\alpha_{ij}=\alpha[i,j],$
\[\alpha[\alpha[i,j],k]=\alpha[\alpha[i,k],\alpha[j,k]].\]} 
\end{list}
\end{lemma}

\begin{proof}
Suppose $Q$ is a quandle and consider the matrix $M_Q=(\alpha_{ij})$. Since 
$x \tr x=x$, we have 
\[\alpha_{ii}= x_i\tr x_i=x_i.\]
Distinctness of elements on the diagonal is then equivalent to distinctness of
elements of the quandle. Conversely, if $x_i$ appears in two positions on the 
diagonal, say $\alpha_{ii}= x_i=\alpha_{jj}$ then we have $x_j\tr x_j \ne x_j$
and $Q$ is not a quandle.

If $Q$ is a quandle, then since $\alpha_{ij}=x_i\tr x_j$, column $j$ of $M_Q$
consists of elements of the form $x_i \tr x_j$. Quandle axiom (ii) says 
that for every $a,b\in Q$ there is a unique $c$ such that $a=c\tr b$, so 
\[\alpha_{ij}=x_i\tr x_j=x_k\tr x_j=\alpha_{kj}\] 
implies $x_i=x_k$, which implies $i=k$. Conversely, if the entries in  
column $c(x_j)$ are distinct, the fact that there are $n$ entries chosen from
$\{x_1,x_2,\dots, x_n\}$ implies that every element appears in the column 
$c(x_j)$, that is, every element is $x_i\tr x_j$ for a unique $x_i$.

Finally, condition (iii) is simply quandle axiom (iii) rewritten with the
notation $\alpha_{r(x_i)x(x_j)}=x_i\tr x_j$.
\end{proof}

\begin{corollary} \label{rcl}
If $M_Q$ is a quandle matrix, we can read the row and column labels 
off the diagonal: if $\alpha_{ii}=x$, then the entries in row $i$
are of the form $x\tr y$ and the entries in column $i$ are of the form 
$y\tr x$.
\end{corollary}


It is worth noting that if $Q$ is not a quandle but a \textit{rack}, i.e.
if $\{Q,\tr\}$ satisfies quandle axioms (ii) and (iii) but not necessarily 
(i), then corollary \ref{rcl} does not hold, and there is no standard form 
matrix presentation for non-quandle racks. In particular, to represent racks
with matrices, we need to keep track of which row and column represent
which element of $Q$, since unlike the quandle case, in a non-quandle rack 
we cannot recover this information from the matrix itself.

\begin{corollary}
If $M=M_Q$ is an integral quandle matrix for a quandle $Q$ then the trace of
$M_Q$ is \[tr(M_Q) = \frac{n(n+1)}{2}.\] 
\end{corollary}

\begin{proof}
By lemma \ref{mqdef}, the diagonal is a permutation of $\{1,2,\dots,n \}$.
Then 
\[ tr(M_Q) = \sum_{x=1}^{n} x = \frac{n(n+1)}{2}.\]
\end{proof}

\begin{definition}
Let $\rho\in \Sigma_n$ be a permutation of $\{1,2,\dots,n\}$. Set
\[\rho(M_Q) = A_{\rho}^{-1}(\rho(\alpha_{ij}))A_{\rho}\]
where $M_Q=(\alpha_{ij})$ and $A_{\rho}$ is the permutation matrix
of $\rho$. Then we say $\rho(M_Q)$ is \textit{p-equivalent} or 
\textit{permutation-equivalent} to $M_Q$, and write $\rho(M_Q)\sim_p M_Q$.
\end{definition}

The fact that $p$-equivalence is an equivalence relation follows from the
fact that $\Sigma_n$ is a group. We now can prove our main theorem.

\begin{theorem} \label{main}
Two integral quandle matrices in standard form determine isomorphic quandles 
iff they are $p$-equivalent by a permutation $\rho\in \Sigma_n$. 
\end{theorem}

\begin{proof}
Let $\rho:Q\to Q'$ be an isomorphism of finite quandles and let $M_Q$, 
$M_{Q'}$ be the standard form integral quandle matrices of $Q$ and $Q'$ 
respectively. Since $\rho$ is a bijection 
$\rho:\{1,2,\dots, n\}\to\{1,2,\dots, n\} $, we have $\rho \in \Sigma_n$.

Then $\rho(i\tr j)=\rho(i)\tr \rho(j)$ says that in the operation table of 
$Q'$, the element in row $r(\rho(i))$ and column $c(\rho(j))$
is $\rho(i\tr j)$; that is, we obtain an operation table for $Q'$ by applying 
the permutation $\rho$ to every element in the table, including the row and 
column labels (which we can recover from the diagonal). Conjugation by the 
permutation matrix of $\rho$ then puts the matrix back in standard form.

Conversely, if $M_{Q'}$ is $p$-equivalent to $M_{Q}$ by a permutation
$\rho$, then the element in row $r(\rho(i))$ and column 
$c(\rho(j))$ in $M_{Q'}$ is $\rho(i\tr j)$, that is, 
\[\rho(i)\tr \rho(j)=\rho(i\tr j)\] and $\rho$ is an isomorphism of quandles.
\end{proof}

\begin{corollary} \label{aut}
The automorphism group of a finite quandle $Q$ of order $n$ is isomorphic to
the subgroup of $\Sigma_n$ which fixes $M_Q$, i.e.,
\[\mathrm{Aut}(Q) \cong \{\rho\in \Sigma_n \ | \ \rho(M_Q)=M_Q\} 
\subseteq \Sigma_n.\]
\end{corollary}

\begin{proof}
A quandle automorphism of $Q$ is a quandle isomorphism $\rho:Q\to Q$. 
Theorem \ref{main} then implies that $\rho\in \Sigma_n$ induces an 
automorphism of $Q$ iff $\rho(M_Q)=M_Q$.
\end{proof}

\begin{example}
The trivial quandle of order $n$, $T_n$, has integral quandle matrix
\[M_{T_n} = 
\left[\begin{array}{cccc} 
1 & 1 & \dots & 1 \\
2 & 2 & \dots & 2 \\ 
\vdots & \vdots & \ddots & \vdots \\
n & n & \dots & n 
\end{array}\right].\]
It is easy to check that $\rho(M_{T_n})=M_{T_n}$ for all $\rho\in \Sigma_n$,
and by corollary \ref{aut}, $\mathrm{Aut}(T_n)\cong \Sigma_n$.
\end{example}

\begin{example}
The quandle matrix 
\[M_Q=\left[\begin{array}{cccc} 
1 & 1 & 1 & 1 \\ 
2 & 2 & 2 & 3 \\ 
3 & 3 & 3 & 2 \\ 
4 & 4 & 4 & 4
\end{array}\right]
\mathrm{ \quad is \ }p\mathrm{-equivalent \ to \quad } \rho(M_Q) = 
\left[\begin{array}{cccc} 
1 & 1 & 2 & 1 \\ 
2 & 2 & 1 & 2 \\ 
3 & 3 & 3 & 3 \\ 
4 & 4 & 4 & 4
\end{array}\right]
\]
with $\rho=(1432)$, since
\begin{eqnarray*}
\left[\begin{array}{cccc} 
1 & 1 & 2 & 1 \\ 
2 & 2 & 1 & 2 \\ 
3 & 3 & 3 & 3 \\ 
4 & 4 & 4 & 4 \end{array}\right] & = & 
\left[ \begin{array}{cccc} 
0 & 1 & 0 & 0 \\ 
0 & 0 & 1 & 0 \\
0 & 0 & 0 & 1 \\ 
1 & 0 & 0 & 0 \\  \end{array}\right]
\left[\begin{array}{cccc} 
4 & 4 & 4 & 4 \\ 
1 & 1 & 1 & 2 \\ 
2 & 2 & 2 & 1 \\ 
3 & 3 & 3 & 3
\end{array}\right]
\left[ \begin{array}{cccc} 
0 & 0 & 0 & 1 \\ 
1 & 0 & 0 & 0 \\
0 & 1 & 0 & 0 \\ 
0 & 0 & 1 & 0 \\  \end{array}\right] \\
& = & A_{\rho}^{-1} (\rho(M_Q)) A_{\rho}.
\end{eqnarray*} 
\end{example}

The number $N_p(Q)$ of standard form integral quandle matrices in the 
$p$-equiv\-alence class of $Q$ is an invariant of quandle type. A 
conjugation $\phi \rho \phi^{-1}$ of an automorphism $\rho \in 
\mathrm{Aut}(Q)$ by an isomorphism $\phi:Q\to Q'$ yields an automorphism 
of $Q'$, so we have 
\[ Q\cong Q' \ \Rightarrow N_p(Q)=N_p(Q').\]
Then since every permutation $\rho\in \Sigma$ defines either an automorphism 
of $Q$ or an isomorphism from $Q$ to a $p$-equivalent quandle $Q'$, we have

\begin{corollary} Let $Q$ be a quandle with $n$ elements. Then
\[|\Sigma_n|=N_p(Q) |\mathrm{Aut}(Q)|.\]
\end{corollary}

Joyce, in \cite{J}, defined quandle to be \textit{algebraically connected} or 
just \textit{connected} if the quandle has only one orbit under the inner 
automorphism group, that is, if the set 
\[O(a) = \{(\dots ((a\diamond_1 b_1) \diamond_2 b_2) \dots \diamond_n b_n) 
\quad | \quad
b_i\in Q, \ \diamond_i\in\{\tr,\tl\}\}=Q \quad \mathrm{for\ all}\  a\in Q.\]

By lemma \ref{mqdef}, we know that the columns in an integral quandle matrix
$M_Q$ must be permutations of $\{1,2,\dots, n\}$. If the rows in $M_Q$ are 
also permutations of $\{1,2,\dots,n\}$, then $Q$ is connected. A 
matrix in which both rows and columns are permutations of $\{1,2,\dots,n\}$ is 
called a \textit{latin square}, and a quandle whose matrix is a latin square
is connected. However, not every latin square is a quandle matrix; 
for example, the latin square
\[\left[\begin{array}{ccc} 
1 & 2 & 3 \\ 
3 & 1 & 2 \\ 
2 & 3 & 1 
\end{array}\right]\]
fails the first condition for being a quandle matrix. 

\begin{definition}
A quandle $Q$ is \textit{latin} if the matrix of $Q$ is a latin square,
that is, if every row of the matrix of $Q$ is a permutation of $\{1,2,\dots, 
n\}$.
\end{definition}

Moreover, not every connected quandle is latin. The conditions of $Q$ being
latin and $Q$ being connected coincide when $Q$ is the conjugation quandle of 
a group, since 
\begin{eqnarray*}
(\dots ((a\diamond_1 b_1) \diamond_2 b_2) \dots \diamond_n b_n) & = &
b_n^{-j_n}\dots (b_2^{-j_2}(b_1^{-j_1}ab_1^{j_1})b_2^{j_2})\dots 
b_n^{j_n} \\ 
& = & (b_n^{-j_n}\dots b_2^{-j_2}b_1^{-j_1})a(b_1^{j_1}b_2^{j_2}\dots 
b_n^{j_n}) \\
& = & (b_1^{j_1}b_2^{j_2} \dots b_n^{j_n})^{-1} a(b_1^{j_1}b_2^{j_2} 
\dots b_n^{j_n}),
\end{eqnarray*}
where $j_k=\pm 1$, and every element of $O(a)$ is $a\tr b$ for some $b\in Q$. 
If a quandle is isomorphic to union of a proper subset of conjugacy classes in 
a group, then the group elements defining some inner automorphisms may not be 
elements of the quandle, and we can have connected quandles which are 
non-latin. For example, the quandle of transpositions in $\Sigma_6$ is 
connected and non-latin.\footnote{Thanks to Steven Wallace for bringing this 
example to the authors' attention.}
\[
M_Q= \left[
\begin{array}{cccccc}
 1 & 4 & 5 & 2 & 3 & 1 \\
 4 & 2 & 6 & 1 & 2 & 3 \\
 5 & 6 & 3 & 3 & 1 & 2 \\
 2 & 1 & 4 & 4 & 6 & 5 \\
 3 & 5 & 1 & 6 & 5 & 4 \\
 6 & 3 & 2 & 5 & 4 & 6
\end{array}
\right]
\]

Connected quandles are of prime interest since knot quandles are connected.
A list of known connected quandles together with an algorithm for finding
connected quandles is given in \cite{O}. A previous computer search by S. 
Yamada for isomorphism classes of quandles is mentioned, though only the 
resulting connected quandles are listed.

\section{Computational results}

In this section, we describe an algorithm for determining all quandles of 
order $n$ by computing all standard form integral quandle matrices of order 
$n$. We then give the results of application of this algorithm for $n=3$, 
$n=4$ and $n=5$. We also determine the automorphism group of each quandle as 
well as a presentation of $Q$ as an Alexander quandle when appropriate.
The maple code used to obtain these results is available on the second
author's website at \texttt{http://www.esotericka.org/quandles}, as is some
more recent and much faster C code \cite{HMN}.

To determine all quandles of order $n$, we first determine for each 
$i=1,\dots, n$ a list $P_{n,i}$ of all vectors whose entries are permutations 
of the set $\{1,2,\dots,n\}$ with entry $i$ in the $i$th position. We then 
consider all matrices $M=(\alpha_{ij})$ with columns $C_1,\dots,C_n$ chosen
from $P_{n,i}$ respectively, since we lose no generality by considering only 
quandle matrices in standard form. For each matrix which satisfies this 
condition, we then check whether 
\[ \alpha_{\alpha_{ij}k}=\alpha_{\alpha_{ik}\alpha_{jk}}\]
for each triple $i,j,k=1,2,\dots, n$.\footnote{It is helpful to make sure 
the program exits the loop at the first triple which does not satisfy the 
condition!} We then check the resulting list of quandle matrices
to determine $p$-equivalence classes. One way to do this is to
compare $M$ and $\rho(M')$ for each $\rho\in P_n$ for every pair $M,M'$ of
quandle matrices, removing $M'$ from the list whenever $M=\rho(M')$ for 
some $M\ne M'$. To compute $\mathrm{Aut}(Q)$, we simply note which 
permutations fix a representative matrix $M_Q$ of $Q$. 

It is easy to check that there is only one quandle of order 1 and one 
quandle of order 2, both trivial (i.e., $x\tr y=x \ \forall x\in Q$.) 
Application of the above algorithm shows that there are three quandle 
isomorphism classes of order 3, 7 isomorphism classes of quandles of order 4 
and 22 isomorphism classes of quandles of order 5. Representative quandle 
matrices for each of these are listed in the tables below. In general, for 
quandles of order $n$, the above algorithm requires $(n-1)!^n$ passes 
through the loop, each pass of which can require up to $n^3$ checks 
of the third quandle condition. 

As a question for further research, we would like to know whether 
there are quandle invariants derivable from $M_Q$ via linear algebra.
A natural first attempt to find such an invariant is to consider the
determinant of $M_Q$. Unfortunately, $p$-equivalence does not generally 
preserve determinants:

\[\left|\begin{array}{ccccc}
1 & 4 & 5 & 2 & 3 \\
3 & 2 & 1 & 5 & 4 \\
4 & 5 & 3 & 1 & 2 \\
5 & 3 & 2 & 4 & 1 \\
2 & 1 & 4 & 3 & 5 \\
\end{array}\right| = -825 \ne -1875 = 
\left|\begin{array}{ccccc}
1 & 5 & 4 & 3 & 2 \\
3 & 2 & 1 & 5 & 4 \\
5 & 4 & 3 & 2 & 1 \\
2 & 1 & 5 & 4 & 3 \\
4 & 3 & 2 & 1 & 5 \\
\end{array}\right|,
\] 
but these two matrices are $p$-equivalent via the permutation $(153)(24)$.

\pagebreak

\begin{figure}[!ht] 
\begin{center}
\begin{tabular}{|cccc|} \hline 
$M_Q$ &  Alexander presentation & $\mathrm{Aut}(Q)$ & $N_p(Q)$
\\ \hline
& & &  \\
$\left[\begin{array}{ccc} 
1 & 1 & 1 \\ 
2 & 2 & 2 \\ 
3 & 3 & 3 
\end{array}\right]$ & $\Lambda_3/(t+2)$ & $\Sigma_3$ & 1 \\
& & &  \\
$\left[\begin{array}{ccc} 
1 & 3 & 2 \\ 
3 & 2 & 1 \\ 
2 & 1 & 3  
\end{array}\right]$ & $\Lambda_3/(t+1)$ & $\Sigma_3$ & 1 \\ 
& & &  \\
$\left[\begin{array}{ccc} 
1 & 1 & 1 \\ 
3 & 2 & 2 \\ 
2 & 3 & 3  
\end{array}\right]$ 
 & -- & $\mathbb{Z}_2$ & 3 \\ 
& & &  \\ \hline
\end{tabular}
\caption{Quandle matrices for quandles of order 3}
\end{center}
\end{figure}

\begin{figure}
\begin{center}
\begin{tabular}{|ccc|} \hline
$M_Q$ & Alexander presentation &  $\mathrm{Aut}(Q)$  \\ \hline
& &   \\
$\left[\begin{array}{cccc} 
1 & 1 & 1 & 1 \\ 
2 & 2 & 2 & 2 \\ 
3 & 3 & 3 & 3 \\ 
4 & 4 & 4 & 4 \\ 
\end{array}\right]$ & $\Lambda_4/(t+3)$ &$\Sigma_4$  \\
& &   \\
$\left[\begin{array}{cccc} 
1 & 1 & 1 & 1 \\ 
2 & 2 & 2 & 3 \\ 
3 & 3 & 3 & 2 \\ 
4 & 4 & 4 & 4 \\ 
\end{array}\right]$ & -- & $\mathbb{Z}_2$   \\
& &  \\
$\left[\begin{array}{cccc} 
1 & 1 & 1 & 2 \\ 
2 & 2 & 2 & 3 \\ 
3 & 3 & 3 & 1 \\ 
4 & 4 & 4 & 4 \\ 
\end{array}\right]$ & -- & $\mathbb{Z}_3$  \\
& &   \\
$\left[\begin{array}{cccc} 
1 & 1 & 2 & 2 \\ 
2 & 2 & 1 & 1 \\ 
3 & 3 & 3 & 3 \\ 
4 & 4 & 4 & 4 \\ 
\end{array}\right]$ & -- & $\mathbb{Z}_2\oplus\mathbb{Z}_2$  \\
& &  \\
$\left[\begin{array}{cccc} 
1 & 1 & 1 & 1 \\ 
2 & 2 & 4 & 3 \\ 
3 & 4 & 3 & 2 \\ 
4 & 3 & 2 & 4 \\ 
\end{array}\right]$ & -- & $\Sigma_3$  \\
& &   \\
$\left[\begin{array}{cccc} 
1 & 1 & 2 & 2 \\ 
2 & 2 & 1 & 1 \\ 
4 & 4 & 3 & 3 \\ 
3 & 3 & 4 & 4 \\ 
\end{array}\right]$ & $\Lambda_2/(t^2+1)$ & $D_8$ \\
& &   \\
$\left[\begin{array}{cccc} 
1 & 4 & 2 & 3 \\ 
3 & 2 & 4 & 1 \\ 
4 & 1 & 3 & 2 \\ 
2 & 3 & 1 & 4 \\ 
\end{array}\right]$ & $\Lambda_2/(t^2+t+1)$ & $A_4$   \\
& &   \\ \hline
\end{tabular}
\caption{Quandle matrices for quandles of order 4}
\end{center}
\end{figure}

\begin{figure} 
\begin{center}
{\small
\begin{tabular}{|ccc|ccc|} \hline
$Q_M$ & $\genfrac{}{}{0pt}{0}{\textrm{Alex.}}{\textrm{pres.}}$ & $\mathrm{Aut}(Q)$ & 
$Q_M$ & $\genfrac{}{}{0pt}{0}{\textrm{Alex.}}{\textrm{pres.}}$ & $\mathrm{Aut}(Q)$ \\  \hline
 & & & & & \\
$\left[\begin{array}{ccccc}
1 & 1 & 1 & 1 & 1 \\
2 & 2 & 2 & 2 & 2 \\
3 & 3 & 3 & 3 & 3 \\
4 & 4 & 4 & 4 & 4 \\
5 & 5 & 5 & 5 & 5 \\
\end{array}\right]$ & \hspace{-3pt}$\Lambda_5/(t+4)$\hspace{-7.2pt} & $\Sigma_5$ & 
$\left[\begin{array}{ccccc}
1 & 1 & 1 & 1 & 1 \\
2 & 2 & 2 & 2 & 2 \\
3 & 3 & 3 & 3 & 4 \\
4 & 4 & 4 & 4 & 3 \\
5 & 5 & 5 & 5 & 5 \\
\end{array}\right]$ & -- & $\mathbb{Z}_2\oplus\mathbb{Z}_2$ \\ 
 & & & & & \\
$\left[\begin{array}{ccccc}
1 & 1 & 1 & 1 & 1 \\
2 & 2 & 2 & 2 & 3 \\
3 & 3 & 3 & 3 & 4 \\
4 & 4 & 4 & 4 & 2 \\
5 & 5 & 5 & 5 & 5 \\
\end{array}\right]$ & -- & $\mathbb{Z}_3$ & 
$\left[\begin{array}{ccccc}
1 & 1 & 1 & 1 & 2 \\
2 & 2 & 2 & 2 & 1 \\
3 & 3 & 3 & 3 & 4 \\
4 & 4 & 4 & 4 & 3 \\
5 & 5 & 5 & 5 & 5 \\
\end{array}\right]$ & -- & $D_8$ \\
 & & & & & \\
$\left[\begin{array}{ccccc}
1 & 1 & 1 & 1 & 2 \\
2 & 2 & 2 & 2 & 3 \\
3 & 3 & 3 & 3 & 4 \\
4 & 4 & 4 & 4 & 1 \\
5 & 5 & 5 & 5 & 5 \\
\end{array}\right]$ & -- & $\mathbb{Z}_4$ & 
$\left[\begin{array}{ccccc}
1 & 1 & 1 & 1 & 1 \\
2 & 2 & 2 & 3 & 3 \\
3 & 3 & 3 & 2 & 2 \\
4 & 4 & 4 & 4 & 4 \\
5 & 5 & 5 & 5 & 5 \\
\end{array}\right]$ & -- & $\mathbb{Z}_2\oplus\mathbb{Z}_2$ \\
 & & & & & \\
$\left[\begin{array}{ccccc}
1 & 1 & 1 & 2 & 2 \\
2 & 2 & 2 & 3 & 3 \\
3 & 3 & 3 & 1 & 1 \\
4 & 4 & 4 & 4 & 4 \\
5 & 5 & 5 & 5 & 5 \\
\end{array}\right]$ & -- & $\mathbb{Z}_3\oplus\mathbb{Z}_2$ & 
$\left[\begin{array}{ccccc}
1 & 1 & 1 & 2 & 3 \\
2 & 2 & 2 & 3 & 1 \\
3 & 3 & 3 & 1 & 2 \\
4 & 4 & 4 & 4 & 4 \\
5 & 5 & 5 & 5 & 5 \\
\end{array}\right]$ & -- & $\Sigma_3$ \\
 & & & & & \\
$\left[\begin{array}{ccccc}
1 & 1 & 1 & 1 & 1 \\
2 & 2 & 2 & 2 & 2 \\
3 & 3 & 3 & 5 & 4 \\
4 & 4 & 5 & 4 & 3 \\
5 & 5 & 4 & 3 & 5 \\
\end{array}\right]$ & -- & $\Sigma_3\times \mathbb{Z}_2$ & 
$\left[\begin{array}{ccccc}
1 & 1 & 1 & 2 & 2 \\
2 & 2 & 2 & 1 & 1 \\
3 & 3 & 3 & 3 & 3 \\
4 & 4 & 5 & 4 & 4 \\
5 & 5 & 4 & 5 & 5 \\
\end{array}\right]$ & -- & $\mathbb{Z}_2\oplus\mathbb{Z}_2$ \\ 
 & & & & & \\
$\left[\begin{array}{ccccc}
1 & 1 & 2 & 2 & 2 \\
2 & 2 & 1 & 1 & 1 \\
3 & 3 & 3 & 3 & 3 \\
4 & 4 & 4 & 4 & 4 \\
5 & 5 & 5 & 5 & 5 \\
\end{array}\right]$ & -- & $\Sigma_3\times \mathbb{Z}_2$ & 
$\left[\begin{array}{ccccc}
1 & 1 & 2 & 2 & 2 \\
2 & 2 & 1 & 1 & 1 \\
3 & 3 & 3 & 3 & 4 \\
4 & 4 & 4 & 4 & 3 \\
5 & 5 & 5 & 5 & 5 \\
\end{array}\right]$ & -- & $\mathbb{Z}_2\oplus\mathbb{Z}_2$ \\ 
 & & & & & \\
\hline
\end{tabular}
}
\caption{Quandle matrices for quandles of order 5 - part 1}
\end{center}
\end{figure}

\begin{figure} 
\begin{center}
{\small
\begin{tabular}{|ccc|ccc|} \hline
$Q_M$ & $\genfrac{}{}{0pt}{0}{\textrm{Alex.}}{\textrm{pres.}}$ & $\mathrm{Aut}(Q)$ & 
$Q_M$ & $\genfrac{}{}{0pt}{0}{\textrm{Alex.}}{\textrm{pres.}}$ & $\mathrm{Aut}(Q)$ \\  \hline
 & & & & & \\
$\left[\begin{array}{ccccc}
1 & 1 & 2 & 2 & 2 \\
2 & 2 & 1 & 1 & 1 \\
3 & 3 & 3 & 5 & 4 \\
4 & 4 & 5 & 4 & 3 \\
5 & 5 & 4 & 3 & 5 \\
\end{array}\right]$ & -- & \hspace{-6pt}$\Sigma_3\times \mathbb{Z}_2$ & 
$\left[\begin{array}{ccccc}
1 & 1 & 2 & 2 & 2 \\
2 & 2 & 1 & 1 & 1 \\
3 & 3 & 3 & 3 & 3 \\
5 & 5 & 5 & 4 & 4 \\
4 & 4 & 4 & 5 & 5 \\
\end{array}\right]$ & -- & $D_8$ \\
 & & & & & \\
$\left[\begin{array}{ccccc}
1 & 1 & 1 & 1 & 1 \\
2 & 2 & 5 & 3 & 4 \\
3 & 4 & 3 & 5 & 2 \\
4 & 5 & 2 & 4 & 3 \\
5 & 3 & 4 & 2 & 5 \\
\end{array}\right]$ & -- & $A_4$ & 
$\left[\begin{array}{ccccc}
1 & 1 & 1 & 1 & 1 \\
2 & 2 & 2 & 3 & 3 \\
3 & 3 & 3 & 2 & 2 \\
5 & 5 & 5 & 4 & 4 \\
4 & 4 & 4 & 5 & 5 \\
\end{array}\right]$ & -- & \hspace{-5pt}$\mathbb{Z}_2\oplus\mathbb{Z}_2$ \\
 & & & & & \\
$\left[\begin{array}{ccccc}
1 & 1 & 1 & 2 & 2 \\
2 & 2 & 2 & 3 & 3 \\
3 & 3 & 3 & 1 & 1 \\
5 & 5 & 5 & 4 & 4 \\
4 & 4 & 4 & 5 & 5 \\
\end{array}\right]$ & -- & \hspace{-6pt}$\mathbb{Z}_3\oplus\mathbb{Z}_2$ & 
$\left[\begin{array}{ccccc}
1 & 3 & 4 & 5 & 2 \\
3 & 2 & 5 & 1 & 4 \\
4 & 5 & 3 & 2 & 1 \\
5 & 1 & 2 & 4 & 3 \\
2 & 4 & 1 & 3 & 5 \\
\end{array}\right]$ & $\Lambda_5/(t+2)$\hspace{-6pt} & $D_{20}$ \\
 & & & & & \\
$\left[\begin{array}{ccccc}
1 & 1 & 2 & 2 & 2 \\
2 & 2 & 1 & 1 & 1 \\
4 & 5 & 3 & 5 & 4 \\
5 & 3 & 5 & 4 & 3 \\
3 & 4 & 4 & 3 & 5 \\
\end{array}\right]$ & -- & $\Sigma_3$ & 
$\left[\begin{array}{ccccc}
1 & 4 & 5 & 3 & 2 \\
3 & 2 & 4 & 5 & 1 \\
2 & 5 & 3 & 1 & 4 \\
5 & 1 & 2 & 4 & 3 \\
4 & 3 & 1 & 2 & 5 \\
\end{array}\right]$ & $\Lambda_5/(t+1)$\hspace{-6pt} & $D_{20}$ \\
 & & & & & \\
$\left[\begin{array}{ccccc}
1 & 1 & 1 & 1 & 1 \\
2 & 2 & 2 & 3 & 3 \\
3 & 3 & 3 & 2 & 2 \\
4 & 5 & 5 & 4 & 4 \\
5 & 4 & 4 & 5 & 5 \\
\end{array}\right]$ & -- & $D_8$ & 
$\left[\begin{array}{ccccc}
1 & 4 & 5 & 2 & 3 \\
3 & 2 & 1 & 5 & 4 \\
4 & 5 & 3 & 1 & 2 \\
5 & 3 & 2 & 4 & 1 \\
2 & 1 & 4 & 3 & 5 \\
\end{array}\right]$ & $\Lambda_5/(t+3)$\hspace{-6pt} & $D_{20}$ \\
 & & & & & \\
\hline
\end{tabular}
}
\caption{Quandle matrices for quandles of order 5 - part 2}
\end{center}
\end{figure}

\pagebreak

\vskip+1cm
\par
\noindent
This article is available at
\texttt{http://intlpress.com/HHA/v7/n1/a11/}
\vskip+0.5cm

\end{document}